\documentclass[12pt,a4paper]{article}
\usepackage{amsmath,amsxtra,amsfonts,amssymb,latexsym,amsthm,amscd,verbatim,amstext}
\usepackage[mathscr]{eucal}

\newcommand{\dif}{\mathrm{d}}

\input cyracc.def

\newfont{\tricyr}{wncyr10 at 12pt}
\newfont{\tricyi}{wncyi10 at 12pt}
\newfont{\tricyb}{wncyb10 at 12pt}
\newfont{\Tricyr}{wncyr10 at 13.6pt}
\newfont{\Tricyi}{wncyi10 at 13.6pt}
\newfont{\Tricyb}{wncyb10 at 13.6pt}
\newfont{\tricmr}{cmr10 at 13.6pt}
\newfont{\tricmi}{cmti10 at 13.6pt}
\newfont{\tricmb}{cmb10 at 13.6pt}

\theoremstyle{plain}
\newtheorem{Th}{Theorem}

\newtheorem{Lem}{Lemma}

\theoremstyle{definition}

\newtheorem{Not}{Remark}

\begin{document}
%%\NoBlackBoxes%%
\centerline {{\bf The existence and decay rates of strong solutions }} 
\centerline {{\bf for Navier-Stokes Equations in Sobolev}} 

\vskip 0.7cm
\begin{center}
D. Q. Khai
\end{center}
\begin{center}
Institute of Mathematics, VAST\\
18 Hoang Quoc Viet, 10307 Cau Giay, Hanoi, Vietnam
\end{center}
\vskip 0.7cm 
{\bf Abstract}: In this paper, we prove some results on the existence 
and decay properties of high order derivatives
in time and space variables for local and global solutions of the Cauchy problem for the
Navier-Stokes equations in Sobolev.
\footnotetext[3]{{\it e-mail address}: khaitoantin@gmail.com, triminh@math.ac.vn}
\vskip 0,5cm 
\centerline{\bf \S1.  Introduction} 
\vskip 0,5cm 
This paper studies the Cauchy problem of the incompressible Navier-Stokes equations (NSE)
in the whole space $\mathbb R^d$ for $d\geq 2$,
\begin{align} 
\left\{\begin{array}{ll} \partial _tu  = 
\Delta u - \nabla .(u \otimes u) - \nabla p , & \\ 
\nabla .u = 0, & \\
u(0, x) = u_0, 
\end{array}\right . \notag
\end{align}
which is a condensed writing for
\begin{align} 
\left\{\begin{array}{ll} 1 \leq k \leq d, \ \  \partial _tu_k  
= \Delta u_k - \sum_{l =1}^{d}\partial_l(u_lu_k) - \partial_kp , & \\ 
\sum_{l =1}^{d}\partial_lu_l = 0, & \\
1 \leq k \leq d, \ \ u_k(0, x) = u_{0k} .
\end{array}\right . \notag
\end{align}
The unknown quantities are the velocity 
$u(t, x)=(u_1(t, x),\dots,u_d(t, x))$ of the fluid 
element at time $t$ and position $x$ and the pressure $p(t, x)$.\\
There is an extensive literature on the existence and decay rate of strong 
\centerline{\bf \S3. Tools from harmonic analysis} 
\vskip 0.5cm
\centerline{\bf \S2. Statement of the results} 
\vskip 0,5cm
Now, for $T > 0$, we say that $u$ is a mild solution of NSE on $[0, T]$ 
corresponding to a divergence-free initial datum $u_0$ 
when $u$ solves the integral equation
$$
u = e^{t\Delta}u_0 - \int_{0}^{t} e^{(t-\tau) \Delta} \mathbb{P} 
\nabla  .\big(u(\tau,.)\otimes u(\tau,.)\big) \dif\tau.
$$
Above we have used the following notation: 
for a tensor $F = (F_{ij})$ 
we define the vector $\nabla.F$ by 
$(\nabla.F)_i = \sum_{j = 1}^d\partial_jF_{ij}$ 
and for two vectors $u$ and $v$, we define their tensor 
product $(u \otimes v)_{ij} = u_iv_j$. 
The operator $\mathbb{P}$ is the Helmholtz-Leray 
projection onto the divergence-free fields 
$$
(\mathbb{P}f)_j =  f_j + \sum_{1 \leq k \leq d} R_jR_kf_k, 
$$
where $R_j$ is the Riesz transforms defined as 
$$
R_j = \frac{\partial_j}{\sqrt{- \Delta}}\ \ {\rm i.e.} \ \  
\widehat{R_jg}(\xi) = \frac{i\xi_j}{|\xi|}\hat{g}(\xi).
$$
The heat kernel $e^{t\Delta}$ is defined as 
$$
e^{t\Delta}u(x) = ((4\pi t)^{-d/2}e^{-|.|^2/4t}*u)(x).
$$
For a space of functions defined on $\mathbb R^d$, say $E(\mathbb R^d)$, we will abbreviate it as $E$ and we do not distinguish between the vector-valued and scalar-value spaces of functions. Throughout the paper, we sometimes use the notation $A \lesssim B$ as an equivalent to $A \leq CB$ with a uniform constant $C$. The notation $A \simeq B$
means that $A \lesssim B$ and $B \lesssim A$. Let $\beta \geq 0$, we define the space $L^\infty(|x|^\beta{\rm dx}):= L^\infty(\mathbb R^d,|x|^\beta{\rm dx})$ which is made up by the measurable  functions $u$ such that  
$$
\big\|u\big\|_{L^\infty(|x|^\beta{\rm dx})}: = \underset{x \in \mathbb R^d}{\rm esssup}|x|^\beta|u(x)| < +\infty.
$$
Now we can state our result
\begin{Th}\label{th4} 
	Assume that $d\geq 1$, and $0\leq \gamma \leq 1 \leq \beta <d$. Then for all $f \in  L^\infty(|x|^\gamma{\rm dx}) \cap L^\infty(|x|^\beta{\rm dx})$ we have
	\begin{gather*}
		\underset{x \in \mathbb R^d, t>0}{\rm sup}\big(|x|^{\tilde \gamma}t^{\frac{1}{2}(\gamma-  \tilde \gamma)} + |x|^\alpha t^{\frac{1}{2}(1-  \alpha)}+|x|^{\tilde\beta} t^{\frac{1}{2}(\beta-  \tilde \beta)}\big)|e^{t\Delta}f| \\ \lesssim  \big\|f\big\|_{L^\infty(|x|^\gamma{\rm dx})}+\big\|f\big\|_{L^\infty(|x|^\beta{\rm dx})}
	\end{gather*}
	for $0 \leq \tilde\gamma \leq \gamma, 0\leq \alpha \leq 1$, and $0 \leq \tilde \beta \leq \beta$.
\end{Th}

\begin{Th}\label{th2}   Let $0\leq \gamma \leq 1\leq \beta <d$ be fixed, then for all $\tilde \gamma, \alpha$, and $\tilde \beta$ satisfying 
	$$
	0\leq\tilde \gamma \leq\gamma, \tilde \beta \geq 0, \beta -2 < \tilde \beta \leq \beta, 0 < \alpha <1, \  and\ \beta -\tilde \beta -1< \alpha < d-\tilde\beta,
	$$ 
	there exists a positive constant $\delta_{\gamma,\tilde \gamma,\alpha,\beta,\tilde \beta,d}$ such that for all
	$u_0 \in L^\infty(|x|^\gamma{\rm d}x) \cap L^\infty(|x|^\beta{\rm d}x)$  with ${\rm div}(u_0) = 0$ satisfying
	\begin{equation}\label{eq4}
		\underset{x \in \mathbb R^d, t>0}{\rm sup}\big(|x|^{\tilde \gamma}t^{\frac{1}{2}(\gamma-  \tilde \gamma)} + |x|^\alpha t^{\frac{1}{2}(1-  \alpha)}+|x|^{\tilde\beta} t^{\frac{1}{2}(\beta-  \tilde \beta)}\big)|e^{t\Delta}u_0|
		\leq \delta_{\gamma,\tilde \gamma,\alpha,\beta,\tilde \beta,d},
	\end{equation}
	NSE has a global mild solution $u$ on $(0,\infty) \times \mathbb R^d$ such that
	\begin{equation}\label{eq5}
		\underset{x \in \mathbb R^d, t>0}{\rm sup}\big(|x|^{\gamma}+t^{\frac{\gamma}{2}} +|x|^\beta +t^{\frac{\beta}{2}}\big)|u(x,t)| <+\infty.
	\end{equation}
\end{Th}
\begin{Not} Our result improves the previous result for $L^\infty(\mathbb R^d, (1+|x|)^\beta{\rm dx})$. This space,
	studied in  \cite{T. Miyakawa}, is a particular case of the 
	space $L^\infty(|x|^\gamma{\rm d}x) \cap L^\infty(|x|^\beta{\rm d}x)$ when $\gamma=0$. 
	Furthermore, we prove that Takahashi's result holds true under a much weaker condition on the initial data. Indeed, from Lemma \ref{lem4} and Theorem \ref{th4}, it is easily seen that the condition \eqref{eq4} of \linebreak Theorem \ref{th2} is weaker than the condition \eqref{eq1}.
\end{Not}
\begin{Not}We invoke Theorem \ref{th4} to deduce that  if $u_0 \in L^\infty(|x|^\gamma{\rm dx}) \cap L^\infty(|x|^\beta{\rm dx})$ and $\big\|u_0\big\|_{L^\infty(|x|^\gamma{\rm dx})}+ \big\|u_0\big\|_{L^\infty(|x|^\beta{\rm dx})}$ is small enough then the \linebreak condition \eqref{eq4} of Theorem \ref{th2} is valid.
\end{Not}
\begin{Th}\label{th2a}   Let $1\leq \beta <d$ be fixed, then for all $\alpha$ satisfying $0<\alpha< 1$, there exists a positive constant $\delta_{\alpha,d}$ such that for all
	$u_0 \in L^\infty(|x|{\rm d}x) \cap L^\infty(|x|^\beta{\rm d}x)$  with ${\rm div}(u_0) = 0$ satisfying
	\begin{equation}\label{eq4a}
		\underset{x \in \mathbb R^d, t>0}{\rm sup}|x|^\alpha t^{\frac{1}{2}(1-  \alpha)}|e^{t\Delta}u_0|
		\leq \delta_{\alpha,d},
	\end{equation}
	NSE has a global mild solution $u$ on $(0,\infty) \times \mathbb R^d$ such that
	$$
	\underset{x \in \mathbb R^d, t>0}{\rm sup}\big(|x|+t^{\frac{1}{2}}\big)|u(x,t)| <+\infty 
	$$
	and
	$$
	\underset{x \in \mathbb R^d, 0<t<T}{\rm sup}|x|^\beta|u(x,t)| <+\infty,\  for\ all\ T \in (0, \infty).
	$$
\end{Th}
\begin{Not}We invoke Theorem \ref{th4} to deduce that  if $u_0 \in L^\infty(|x|{\rm dx})$ and $\big\|u_0\big\|_{L^\infty(|x|{\rm dx})}$ is small enough then the condition \eqref{eq4a} of Theorem \ref{th2a} is valid.
\end{Not}
\vskip 0.5cm
\centerline{\bf \S3. Some auxiliary results} 
\vskip 0.5cm
In this section we establish some auxiliary lemmas. 
We first prove a version of Young's inequality type 
for convolutions in $L^\infty(|x|^\beta{\rm dx})$ spaces. 
\begin{Lem}\label{lem1} 
	Assume that $d\geq 1, 0<\alpha<d, 0<\beta<d$ and $\alpha +\beta >d$. 
	Then for all  $f \in L^\infty(|x|^\alpha{\rm dx})$ 
	and for all  $g \in L^\infty(|x|^\beta{\rm dx})$ 
	we have
	$$
	\big\|f*g\big\|_{L^\infty(|x|^{\alpha+\beta-d}{\rm dx})} 
	\lesssim \big\|f\big\|_{L^\infty(|x|^{\alpha}{\rm dx})}
	\big\|g\big\|_{L^\infty(|x|^{\beta}{\rm dx})}.
	$$
\end{Lem}
\textbf{Proof}.  Since $f*g$ is bilinear on 
$L^\infty(|x|^\alpha{\rm dx}) \times L^\infty(|x|^\beta{\rm dx})$, 
we may assume $\big\|f\big\|_{L^\infty(|x|^\alpha{\rm dx})} 
=\big\|g\big\|_{L^\infty(|x|^\beta{\rm dx})}=1$. We have
\begin{gather*}
	(f*g)(x) =\int_{\mathbb R^d}f(x-y)g(y){\rm d}y 
	= \int_{|y|<\frac{|x|}{2}} + \int_{\frac{|x|}{2}<|y|<\frac{3|x|}{2}} +\int_{|y|>\frac{3|x|}{2}}
	=I_1+I_2+I_3.
\end{gather*}
From 
$$
|f(x)| \leq |x|^{-\alpha},\ {\rm and}\ |g(x)| \leq |x|^{-\beta},
$$
we get
\begin{gather*}
	|I_1| \leq \int_{|y|<\frac{|x|}{2}}\frac{{\rm d}y}{|x-y|^\alpha |y|^\beta} \leq
	\frac{2^\alpha}{|x|^\alpha}\int_{|y|<\frac{|x|}{2}}\frac{{\rm d}y}{ |y|^\beta}\simeq 
	\frac{1}{|x|^{\alpha+\beta-d}}.
\end{gather*}
\begin{gather*}
	|I_2| \leq \int_{\frac{|x|}{2}<|y|<\frac{3|x|}{2}}\frac{{\rm d}y}{|x-y|^\alpha |y|^\beta} \leq
	\frac{2^\beta}{|x|^\beta}\int_{|y|<\frac{5|x|}{2}}\frac{{\rm d}y}{ |y|^\alpha}\simeq 
	\frac{1}{|x|^{\alpha+\beta-d}}.
\end{gather*}
\begin{gather*}
	|I_3| \leq \int_{|y|>\frac{3|x|}{2}}\frac{{\rm d}y}{|x-y|^\alpha |y|^\beta} \leq 
	3^\alpha\int_{|y|>\frac{3|x|}{2}}\frac{{\rm d}y}{|y|^\alpha |y|^\beta}\simeq 
	\frac{1}{|x|^{\alpha+\beta-d}}.
\end{gather*}
We thus obtain
$$
|(f*g)(x)| \lesssim \frac{1}{|x|^{\alpha+\beta-d}}.
$$
The proof Lemma \ref{lem1} is complete. \qed\\
We now deduce the $L^\infty(|x|^\gamma{\rm dx})-L^\infty(|x|^\beta{\rm dx})$ 
estimate for the heat semigroup.
\begin{Lem}\label{lem2} 
	Assume that $d\geq 1$ and  $0\leq \gamma \leq \beta <d$. 
	Then for all \linebreak $f \in L^\infty(|x|^\beta{\rm dx})$ we have
	\begin{equation}\label{eq6}
		\big\|e^{t\Delta}f\big\|_{L^\infty(|x|^\gamma{\rm dx})} 
		\lesssim t^{-\frac{1}{2}(\beta-\gamma)}
		\big\|f\big\|_{L^\infty(|x|^\beta{\rm dx})},\ for\ t>0.
	\end{equation}
\end{Lem}
\textbf{Proof}. We have 
$$
(e^{t\Delta}f)(x)=\int_{\mathbb R^d}\frac{1}{t^{d/2}}
E\big(\frac{x-y}{\sqrt{t}}\big)f(y){\rm d}y,\ {\rm where}\ E(x) 
= (4\pi )^{-\frac{d}{2}}e^{-\frac{|x|^2}{4}}. 
$$
Recall the simate
\begin{equation}\label{eq7}
	t^{-\frac{d}{2}}e^{-\frac{|x|^2}{4t}} 
	\lesssim |x|^{-\alpha}t^{-\frac{1}{2}(d-\alpha)}, \ {\rm for}\  0 \leq \alpha \leq d.
\end{equation}
We first consider the case $0<\gamma <\beta$. 
From the inequality \eqref{eq7} and Lemma \ref{lem1}, we have
$$
|(e^{t\Delta}f)(x)| \lesssim \int_{\mathbb R^d}
\frac{\big\|f\big\|_{L^\infty(|x|^\beta{\rm dx})}}
{t^{\frac{1}{2}(\beta-\gamma)}|x-y|^{\gamma+d-\beta}
	|y|^{\beta}}{\rm d}y \lesssim 
t^{-\frac{1}{2}(\beta-\gamma)}|x|^{-\gamma}
\big\|f\big\|_{L^\infty(|x|^\beta{\rm dx})}.
$$
This proves \eqref{eq6}.\\
We consider the case $0=\gamma<\beta$. 
Applying Proposition 2.4 $(b)$ in (\cite{P. G. Lemarie-Rieusset 2002}, pp. 20) 
and note that
$|x|^{-\beta} \in L^{\frac{d}{\beta},\infty}$
$$
|e^{t\Delta}f (x)| \lesssim t^{-\frac{d}{2}}\big\|E\big(\frac{.}
{\sqrt t}\big)\big\|_{L^{\frac{d}{d-\beta},1}}\big\|f\big\|_{L^{\frac{d}{\beta},\infty}} \lesssim
t^{-\frac{\beta}{2}}\big\|E\big\|_{L^{\frac{d}{d-\beta},1}}\big\|f\big\|_{L^\infty(|x|^\beta{\rm dx})}.
$$
This proves \eqref{eq6}.\\
Suppose finally that $0 \leq \gamma =\beta$. We have
$$
\int_{\mathbb R^d}\frac{1}{t^{d/2}}E\big(\frac{x-y}
{\sqrt{t}}\big)f(y){\rm d}y = \int_{|y|<\frac{|x|}{2}} 
+ \int_{|y|>\frac{|x|}{2}} = I_1 +I_2.
$$
From the inequality \eqref{eq7}, we have
\begin{gather*}
	|I_1| \lesssim \big\|f\big\|_{L^\infty(|x|^\beta{\rm dx})} 
	\int_{|y|<\frac{|x|}{2}}|x-y|^{-d}|y|^{-\beta}{\rm d}y  \leq \notag \\
	\big\|f\big\|_{L^\infty(|x|^\beta{\rm dx})}\big(\frac{|x|}{2}\big)^{-d}
	\int_{|y|<\frac{|x|}{2}}|y|^{-\beta}{\rm d}y 
	\simeq  \big\|f\big\|_{L^\infty(|x|^\beta{\rm dx})} |x|^{-\beta}.
\end{gather*}
\begin{gather*}
	|I_2| \leq  \big\|f\big\|_{L^\infty(|x|^\beta{\rm dx})}\int_{|y|>\frac{|x|}{2}}
	\frac{1}{t^{d/2}}E\big(\frac{x-y}{\sqrt{t}}\big)|y|^{-\beta}{\rm d}y  \leq \notag\\
	\big\|f\big\|_{L^\infty(|x|^\beta{\rm dx})}\big(\frac{|x|}{2}\big)^{-\beta}
	\int_{y \in \mathbb R^d}\frac{1}{t^{d/2}}E\big(\frac{y}{\sqrt{t}}\big){\rm d}y= 
	C \big\|f\big\|_{L^\infty(|x|^\beta{\rm dx})}|x|^{-\beta}, \notag
\end{gather*}
where
$$
C =2^\beta\int_{y \in \mathbb R^d}E(y){\rm d}y<+\infty.
$$
Therefore,
$$
|e^{t\Delta}f (x)| \lesssim  \big\|f\big\|_{L^\infty(|x|^\beta{\rm dx})}|x|^{-\beta}.
$$
The proof of Lemma \ref{lem2} is complete. \qed \\

We now deduce the $L^\infty(|x|^\gamma{\rm dx})-L^\infty(|x|^\beta{\rm dx})$ 
estimate for the operator $e^{t\Delta} \mathbb{P}\nabla$.
As shown in \cite{P. G. Lemarie-Rieusset 2002}, the kernel function $F_t$ of  $e^{t\Delta} \mathbb{P}\nabla$ satisfies the following \linebreak
inequalities
\begin{equation}\label{eq8}
	F_t(x) = t^{-\frac{d+1}{2}}F\big(\frac{x}{\sqrt{t}}\big), |F(x)| \lesssim \frac{1}{(1+|x|)^{d+1}},
\end{equation}
\begin{equation}\label{eq9}
	|F_t(x)| \lesssim |x|^{-\alpha}t^{-\frac{1}{2}(d+1-\alpha)}, \ {\rm for}\  0 \leq \alpha \leq d+1.
\end{equation}
By using the inequalities \eqref{eq8} and \eqref{eq9} and 
arguing as in the proof of \linebreak Lemma \ref{lem2},  
we can easily prove the following lemma.
\begin{Lem}\label{lem3} 
	Assume that $d\geq 1$ and  $0\leq \gamma \leq \beta <d$. Then for all \linebreak $f \in L^\infty(|x|^\beta{\rm dx})$ we have
	$$
	\big\|e^{t\Delta} \mathbb{P}\nabla.f\big\|_{L^\infty(|x|^\gamma{\rm dx})} \lesssim t^{-\frac{1}{2}(\beta+1-\gamma)}\big\|f\big\|_{L^\infty(|x|^\beta{\rm dx})},\ {\rm for}\ t>0.
	$$
\end{Lem}
\begin{Lem}\label{lem4}
	Let $0\leq \gamma < \beta \leq d$. Assume that $f \in \mathcal{S}'(\mathbb R^d)$ and satisfies the following inequality
	\begin{equation}\label{eq10}
		\underset{x \in \mathbb R^d, t>0}{\rm sup}\big(|x|^{\gamma}+|x|^{\beta}\big)|(e^{t\Delta}f)(x)|  =C <+\infty,
	\end{equation}
	then
	$$
	f \in L^\infty(|x|^\gamma{\rm dx})\cap L^\infty(|x|^\beta{\rm dx})
	$$
	and
	\begin{equation}\label{eq11}
		\underset{x \in \mathbb R^d}{\rm esssup}\big(|x|^{\gamma}+|x|^{\beta}\big)|f(x)|  \leq C.
	\end{equation}
\end{Lem}
\textbf{Proof}. Since $\frac{1}{|x|^{\gamma}+|x|^{\beta}} \in L^{\frac{d}{\beta},\infty}\cap  L^{\frac{d}{\gamma},\infty}$ and $L^{\frac{d}{\beta},\infty}\cap  L^{\frac{d}{\gamma},\infty}\subset L^q$ for all $q$ satisfying $\frac{d}{\beta} < q < \frac{d}{\gamma}$, it follows that $e^{t\Delta}f \in L^{\infty}(0,\infty;L^q)$ for all $q \in \big(\frac{d}{\beta}, \frac{d}{\gamma}\big)$, by a compactness theorem in Banach space, there exists  a sequence $t_k$ which converges to $0$ 
such that $e^{t_k\Delta}f $  converges weakly to $f'$ in $L^q$ with $f' \in L^q$. Since $e^{t\Delta}$ 
is a continuous semigroup on $\mathcal{S}'(\mathbb R^d)$, it follows that $f=f' \in L^q$. Since $e^{t\Delta}$ 
is a continuous semigroup on $L^q(\mathbb R^d), (1\leq q <\infty)$, we get
$$
\underset{k \rightarrow \infty}{\rm lim}\big\|e^{t_k\Delta}f -f\big\|_{L^q}=0,\ {\rm for}\ q\in \big(\frac{d}{\beta}, \frac{d}{\gamma}\big).
$$
Therefore, there exists a subsequence $t_{k_j}$ of the sequence $t_k$ such that
\begin{equation}\label{eq12}
	\underset{j \rightarrow \infty}{\rm lim}(e^{t_{k_j}\Delta}f)(x)=f(x)\ {\rm for\ almost\ everywhere}\ x \in \mathbb{R}^d.
\end{equation}
The inequality \eqref{eq11} is deduced from equalities \eqref{eq10} and \eqref{eq12}. \qed 
\begin{Not} We invoke Lemma \ref{lem4} for $\gamma =0$ and Lemma \ref{lem2} for $\gamma =\beta$  to deduce that the condition \eqref{eq1} of Takahashi on the initial data  
	is equivalent to the condition 
	$$
	\big\|u_0\big\|_{L^\infty((1+|x|)^\beta{\rm dx})} \leq \delta.
	$$
\end{Not}
\begin{Lem}\label{lem5} Let $\gamma, \theta \in \mathbb R$ and $t>0$, then\\
	{\rm (a)} If $\theta< 1$ then
	$$
	\int^{\frac{t}{2}}_0(t-\tau)^{-\gamma} \tau^{-\theta}{\rm d}\tau = C t^{1-\gamma - \theta},\  where\  
	C = \int^{\frac{1}{2}}_0(1-\tau)^{-\gamma} \tau^{-\theta}{\rm d}\tau < \infty.
	$$
	{\rm (b)} If $\gamma < 1$ then
	$$
	\int^{t}_{\frac{t}{2}}(t-\tau)^{-\gamma} \tau^{-\theta}{\rm d}\tau = C t^{1-\gamma - \theta},\  where\  
	C = \int^1_{\frac{1}{2}}(1-\tau)^{-\gamma} \tau^{-\theta}{\rm d}\tau < \infty.
	$$
	{\rm (c)} If $\gamma < 1$ and  $\theta < 1$ then
	$$
	\int^{t}_0(t-\tau)^{-\gamma} \tau^{-\theta}{\rm d}\tau = C t^{1-\gamma - \theta},\  where\  
	C = \int^1_0(1-\tau)^{-\gamma} \tau^{-\theta}{\rm d}\tau < \infty.
	$$
\end{Lem}
The proof of this lemma is elementary and may be omitted.\qed \\
Let us recall the following result on solutions of a quadratic
equation in Banach spaces (Theorem 22.4 in 
\cite{P. G. Lemarie-Rieusset 2002}, p. 227).
\begin{Th}\label{th3}
	Let $E$ be a Banach space, and $B: E \times E \rightarrow  E$ 
	be a continuous bilinear map such that there exists $\eta > 0$ so that
	$$
	\|B(x, y)\| \leq \eta \|x\| \|y\|,
	$$
	for all x and y in $E$. Then for any fixed $y \in E$ 
	such that $\|y\| \leq \frac{1}{4\eta}$, the equation $x = y - B(x,x)$ 
	has a unique solution  $\overline{x} \in E$ satisfying 
	$\|\overline{x}\| \leq \frac{1}{2\eta}$.
\end{Th}
\vskip 0.5cm
\centerline{\bf\S4. Proofs of Theorems \ref{th4}, \ref{th2}, and \ref{th2a}}
\vskip 0.5cm
In this section we will give the proofs of  
Theorems  \ref{th4},  \ref{th2}, and \ref{th2a}. We now need eight more lemmas. 
In order to proceed, we define an auxiliary space $K^\beta_{\alpha,T}$. 
Let $\alpha, \beta$, and $T$ be such that $0 \leq \alpha \leq \beta<d, 0<T\leq +\infty$, 
we define the auxiliary space $K^\beta_{\alpha,T}$ 
which is made up by the measurable functions $u(t,x)$ such that 
\begin{gather*}
	\underset{x \in \mathbb R^d, 0<t<T}{\rm esssup}
	|x|^{\alpha}t^{\frac{1}{2}(\beta-\alpha)}|u(x,t)| < +\infty.
\end{gather*}
The auxiliary space $K^\beta_{\alpha,T}$ is equipped with the norm
$$
\big\|u\big\|_{K^\beta_{\alpha,T}}:= 
\underset{x \in \mathbb R^d, 0<t<T}{\rm esssup}
|x|^{\alpha}t^{\frac{1}{2}(\beta-\alpha)}|u(x,t)|.
$$
We rewrite Lemma  \ref{lem2}  as follows
\begin{Lem}\label{lem6} 
	Assume that $d\geq 1$ and  $0\leq \alpha \leq \beta <d$. 
	Then for all \linebreak  $f \in L^\infty(|x|^\beta{\rm dx})$ we have 
	$e^{t\Delta}f \in K^\beta_{\alpha,T}$ and 
	$\big\|e^{t\Delta}f\big\|_{K^\beta_{\alpha,T}} 
	\leq C\big\|f\big\|_{L^\infty(|x|^\beta{\rm dx})}$, where C is a positive constant independent of T.
\end{Lem}
\begin{Lem}\label{lem7} 
	Assume that $d\geq 1$ and  $0\leq \alpha \leq \beta <d$. 
	Then 
	$$K^\beta_{\alpha,T} \subset K^\beta_{\beta,T} \cap K^\beta_{0,T}.$$
\end{Lem}
The proof of this lemma is elementary and may be omitted.\qed \\
\begin{Lem}\label{lem7a} 
	Assume that $d\geq 1,T<+\infty$, and  $0\leq \alpha \leq \beta \leq \tilde\beta <d$. 
	Then $K^\beta_{\alpha,T}\subset  K^{\tilde\beta}_{\alpha,T}$.
\end{Lem}
The proof of this lemma is elementary and may be omitted.\qed \\
In the following lemmas a particular attention 
will be devoted to the study of the bilinear operator 
$B(u, v)(t)$ defined by 
\begin{equation}\label{eq13}
	B(u, v)(t) = \int_{0}^{t} e^{(t-\tau ) \Delta} \mathbb{P} 
	\nabla.\big(u(\tau)\otimes v(\tau)\big) \dif\tau.
\end{equation}
\begin{Lem}\label{lem8}
	Let $\beta, \tilde \beta, \hat \beta$, and $\alpha$ be such that
	\begin{gather*}
		0 \leq \beta <d, \tilde \beta > \beta -2, 0\leq \tilde \beta \leq \beta, 0<\alpha<1, \beta -\tilde \beta-1 <\alpha < d-\tilde \beta ,\\
		0 \leq \hat \beta \leq \beta,\  and\ \alpha +\tilde \beta-1 < \hat \beta \leq \alpha +\tilde \beta.
	\end{gather*}
	Then the bilinear operator $B$ is continuous from $K^1_{\alpha,T} \times K^\beta_{\tilde \beta,T}$ into $K^\beta_{\hat \beta,T}$
	and the following inequality holds 
	\begin{equation}\label{eq14}
		\big\|B(u, v)\big\|_{K^\beta_{\hat \beta,T}} \leq C\big\|u\big\|_{K^1_{\alpha,T}}
		\big\|v\big\|_{K^\beta_{\tilde \beta,T}},
	\end{equation}
	where C is a positive constant independent of T.
\end{Lem}
\textbf{Proof}. Since $B(.,.)$ is bilinear on 
$K^1_{\alpha,T} \times K^\beta_{\tilde \beta,T}$, 
we may assume  $\big\|u\big\|_{K^1_{\alpha,T}} =\big\|v\big\|_{K^\beta_{\tilde \beta,T}}=1$.
From 
$$
|(u\otimes v)(\tau)| \leq |y|^{-(\alpha+\tilde \beta)}\tau^{-\frac{1}{2}(1-\alpha+\beta-\tilde \beta)}, 
$$
by using Lemma \ref{lem3}, we have
$$
\big|e^{(t-\tau ) \Delta} \mathbb{P}\nabla.\big(u\otimes v\big)\big| 
\lesssim |x|^{-\hat \beta}\frac{1}{(t-\tau)^{\frac{1}{2}
		(1+\alpha+\tilde \beta-\hat \beta)}\tau^{\frac{1}{2}(1-\alpha+\beta-\tilde \beta)}}
$$
then applying Lemma \ref{lem5} (c), we get
$$
|B(u,v)| \lesssim |x|^{-\hat \beta}\int^t_0\frac{1}{(t-\tau)^{\frac{1}{2}
		(1+\alpha+\tilde \beta-\hat \beta)}\tau^{\frac{1}{2}(1-\alpha+\beta-\tilde \beta)}}
{\rm d}\tau \simeq |x|^{-\hat \beta}t^{-\frac{1}{2}(\beta-\hat \beta)}.
$$
This proves Lemma \ref{lem8}.\\
Note that since $\alpha > \beta -\tilde \beta -1$ and $\hat \beta > \alpha+\tilde \beta-1$, 
it follows that the conditions $\frac{1-\alpha +\beta -\tilde \beta}{2} < 1$ and 
$\frac{1+\alpha +\tilde \beta -\hat \beta}{2} < 1$ are valid. So we can apply Lemma \ref{lem5} (c). \qed
\begin{Lem}\label{lem9} Assume that NSE has a mild solution $u \in K^1_{\tilde \alpha,T}$ for some \linebreak $\tilde \alpha \in (0,1)$ 
	with initial data $u_0 \in L^\infty({|x|{\rm d}x})$ then $u \in K^1_{\alpha,T}$ for all $\alpha \in [0,1]$.
\end{Lem}
\textbf{Proof}. From $u = e^{t\Delta}u_0 - B(u,u)$, applying Lemmas \ref{lem6} and \ref{lem8} 
with $\beta =1$ and $\alpha =\tilde \beta =\tilde \alpha$, we get $u \in K^1_{\hat \beta,T}$
for all $\hat \beta \in \big(\tilde \alpha-(1-\tilde \alpha), 2\tilde \alpha\big) \cap[0,1]$.
Applying again Lemmas \ref{lem6} and \ref{lem8}  with $\beta =1, \alpha =\tilde \alpha$, and $\tilde \beta \in \big(\tilde \alpha-(1-\tilde \alpha), 2\tilde \alpha\big) \cap[0,1]$
to get $u \in K^1_{\hat \beta,T}$ for all $\hat \beta \in \big(\tilde \alpha-2(1-\tilde \alpha), 3\tilde \alpha\big) \cap[0,1]$.
By induction, we get $u\in K^1_{\hat \beta,T}$ for all $\hat \beta \in \big(\tilde \alpha-n(1-\tilde \alpha), (n+1)\tilde \alpha\big) \cap[0,1]$ with $n\in \mathbb N$. 
Since $\tilde\alpha \in (0,1)$, it follows that there exists sufficiently large $n$ satisfying
$$
\big(\tilde \alpha-n(1-\tilde \alpha), (n+1)\tilde \alpha\big) \supset [0,1].
$$
This proves Lemma \ref{lem9}. \qed
\begin{Lem}\label{lem10}Let $\beta$ be a fixed number in the interval $[0, d)$. Assume that NSE has a mild solution
	$u \in \underset{\alpha \in [0,1]}{\cap}K^1_{\alpha,T}\cap K^{\beta}_{\tilde \beta,T}$ for some $\tilde \beta \in [0,\beta] \cap (\beta -2, \beta]$ 
	with initial data $u_0 \in L^\infty({|x|^\beta{\rm d}x})$, then $u \in K^\beta_{\hat \beta,T}$ for all $\hat\beta \in [0,\beta]\cap(\tilde \beta-1,\beta]$.
\end{Lem}
\textbf{Proof}. We first prove that $u\in K^\beta_{\hat \beta,T}$ for all $\hat \beta \in [0,\beta]\cap (\tilde \beta -1,\tilde \beta +1)$.\\
Let $\alpha_1$ and $\alpha_2$ be such that
$$
{\rm max}\{\beta -\tilde \beta-1,\hat \beta -\tilde \beta,0\} <\alpha_1<1
$$
and
$$
{\rm max}\{\hat \beta -\tilde \beta,0\} <\alpha_2<{\rm min}\{1,\hat\beta -\tilde\beta +1\}.
$$
We split the integral given in \eqref{eq13} into two parts coming from the subintervals $(0, \frac{t}{2})$ and $(\frac{t}{2}, t)$
\begin{gather*}
	B(u, u)(t) = 
	\int_{0}^{\frac{t}{2}} e^{(t-\tau ) \Delta} \mathbb{P} 
	\nabla.\big(u\otimes u\big) \dif\tau + \int_{\frac{t}{2}}^{t} e^{(t-\tau ) \Delta} \mathbb{P} 
	\nabla.\big(u\otimes u\big) \dif\tau =I_1+I_2.
\end{gather*}
Since $u \in \underset{\alpha \in [0,1]}{\cap}K^1_{\alpha,T}$, it follows that
\begin{equation}\label{eq15}
	|u(x,\tau)| \lesssim |x|^{-\alpha_1}\tau^{-\frac{1}{2}(1-\alpha_1)},
\end{equation}
\begin{equation}\label{eq16}
	|u(x,\tau)| \lesssim |x|^{-\alpha_2}\tau^{-\frac{1}{2}(1-\alpha_2)},
\end{equation}
and since $u\in K^{\beta}_{\tilde \beta,T}$, it follows that
\begin{equation}\label{eq17}
	|u(x,\tau)| \lesssim |x|^{-\tilde \beta}\tau^{-\frac{1}{2}(\beta-\tilde \beta)}.
\end{equation}
From the inequalities \eqref{eq15} and \eqref{eq17}, and Lemma \ref{lem3}, we get
$$
\big|e^{(t-\tau ) \Delta} \mathbb{P}\nabla.\big(u\otimes u\big)\big| \lesssim 
|x|^{-\hat \beta}\frac{1}{(t-\tau)^{\frac{1}{2}(1+\alpha_1+\tilde \beta-\hat \beta)}\tau^{\frac{1}{2}(1-\alpha_1+\beta-\tilde \beta)}}.
$$
Then applying Lemma \ref{lem5} (a), we have
\begin{equation}\label{eq18}
	|I_1| \lesssim |x|^{-\hat \beta}\int^{\frac{t}{2}}_0\frac{1}{(t-\tau)^{\frac{1}{2}(1+\alpha_1+\tilde \beta-\hat \beta)}\tau^{\frac{1}{2}(1-\alpha_1+\beta-\tilde \beta)}}{\rm d}\tau \simeq |x|^{-\hat \beta}t^{-\frac{1}{2}(\beta-\hat \beta)}.
\end{equation}
From the inequalities \eqref{eq16} and \eqref{eq17}, and Lemma \ref{lem3}, we get
$$
\big|e^{(t-\tau ) \Delta} \mathbb{P}\nabla.\big(u\otimes u\big)\big| \lesssim 
|x|^{-\hat \beta}\frac{1}{(t-\tau)^{\frac{1}{2}(1+\alpha_2
		+\tilde \beta-\hat \beta)}\tau^{\frac{1}{2}(1-\alpha_2+\beta-\tilde \beta)}}.
$$
Then applying Lemma \ref{lem5} (b), we have
\begin{equation}\label{eq19}
	|I_2| \lesssim |x|^{-\hat \beta}\int^{t}_{\frac{t}{2}}\frac{1}{(t-\tau)^{\frac{1}{2}
			(1+\alpha_2+\tilde \beta-\hat \beta)}\tau^{\frac{1}{2}
			(1-\alpha_2+\beta-\tilde \beta)}}{\rm d}\tau 
	\simeq |x|^{-\hat \beta}t^{-\frac{1}{2}(\beta-\hat \beta)}.
\end{equation}
From the inequalities \eqref{eq18} and \eqref{eq19}, 
we get $B(u,u)\in K^\beta_{\hat \beta,T}$, 
and from \linebreak $u =e^{t\Delta}u_0+B(u,u)$ and Lemma \ref{lem6}, 
we have $u\in K^\beta_{\hat \beta,T}$.  This proves the result.\\
We now prove  $u \in K^\beta_{\hat \beta,T}$ for all 
$\hat\beta \in [0,\beta] \cap(\tilde\beta-1,\beta]$. Indeed,  
if $\tilde \beta > \beta -1$ then $u\in K^\beta_{\hat \beta,T}$ for all 
$\hat \beta \in [0,\beta]\cap (\tilde \beta -1,\tilde\beta+1) = [0,\beta]\cap (\tilde \beta -1,\beta]$ and so the 
lemma is proved. In the case $\tilde \beta \leq \beta -1$,  in exactly the same way, since $u\in K^\beta_{\hat \beta,T}$ 
for all $\hat \beta \in [0,\beta]\cap (\tilde \beta -1,\tilde\beta+1)$,
it follows that $u\in K^\beta_{\hat \beta,T}$ for all 
$\hat \beta \in [0,\beta]\cap (\tilde \beta -1,\tilde \beta+2)
= [0,\beta]\cap (\tilde \beta -1,\beta]$. Therefore the proof 
of Lemma \ref{lem10} is complete. \qed
\begin{Lem}\label{lem11}Assume that NSE has a mild solution
	$u \in \underset{\alpha \in [0,1]}{\cap}K^1_{\alpha,T}\cap \underset{\hat \beta\in [\tilde\beta,\beta]}{\cap}K^{\beta}_{\hat \beta,T}$ 
	for some $\tilde \beta \in [0,\beta]$ with initial data 
	$u_0 \in L^\infty({|x|^\beta{\rm d}x})$. Then 
	$u \in K^\beta_{\hat \beta,T}$ for all $\hat\beta \in [0,\beta]$.
\end{Lem}
\textbf{Proof}. We first prove that $u\in K^\beta_{\hat \beta,T}$ 
for all $\hat \beta \in [0,\beta]\cap(\tilde \beta-1,\tilde \beta]$.\\
We split the integral given in \eqref{eq13} into two parts coming from the subintervals $(0, \frac{t}{2})$ and $(\frac{t}{2}, t)$
\begin{gather*}
	B(u, u)(t) = 
	\int_{0}^{\frac{t}{2}} e^{(t-\tau ) \Delta} \mathbb{P} 
	\nabla.\big(u\otimes u\big) \dif\tau + \int_{\frac{t}{2}}^{t} e^{(t-\tau ) \Delta} \mathbb{P} 
	\nabla.\big(u\otimes u\big) \dif\tau =I_1+I_2.
\end{gather*}
Let $\alpha_1$ be such that $0<\alpha_1<1$. Since 
$u \in K^1_{\alpha_1,T}\cap K^{\beta}_{\beta,T}$, it follows that
\begin{equation}\label{eq20}
	|u(x,\tau)| \lesssim |x|^{-\alpha_1}\tau^{-\frac{1}{2}(1-\alpha_1)},
\end{equation}
\begin{equation}\label{eq21}
	|u(x,\tau)| \lesssim |x|^{-\beta}.
\end{equation}
From the inequalities \eqref{eq20} and \eqref{eq21}, 
and Lemma \ref{lem3}, we get
$$
\big|e^{(t-\tau ) \Delta} \mathbb{P}\nabla.\big(u\otimes u\big)\big| \lesssim 
|x|^{-\hat \beta}\frac{1}{(t-\tau)^{\frac{1}{2}
		(1+\alpha_1+\beta-\hat \beta)}\tau^{\frac{1}{2}(1-\alpha_1)}}.
$$
Then applying Lemma \ref{lem5} (a), we have
\begin{equation}\label{eq22}
	|I_1| \lesssim |x|^{-\hat \beta}\int^{\frac{t}{2}}_0
	\frac{1}{(t-\tau)^{\frac{1}{2}(1+\alpha_1+\beta-\hat \beta)}
		\tau^{\frac{1}{2}(1-\alpha_1)}}{\rm d}\tau \simeq |x|^{-\hat \beta}
	t^{-\frac{1}{2}(\beta-\hat \beta)}.
\end{equation}
Since $u \in K^1_{0,T}\cap K^{\beta}_{\tilde\beta,T}$, it follows that
\begin{equation}\label{eq23}
	|u(x,\tau)| \lesssim \tau^{-\frac{1}{2}}\ {\rm and}\  |u(x,\tau)| 
	\lesssim |x|^{-\tilde \beta}\tau^{-\frac{1}{2}(\beta-\tilde\beta)}.
\end{equation}
From the inequality \eqref{eq23}, and Lemma \ref{lem3}, we get
$$
\big|e^{(t-\tau ) \Delta} \mathbb{P}\nabla.\big(u\otimes u\big)\big| \lesssim 
|x|^{-\hat \beta}\frac{1}{(t-\tau)^{\frac{1}{2}(1+\tilde\beta-\hat \beta)}
	\tau^{\frac{1}{2}(1+\beta-\tilde\beta)}}.
$$
Then applying Lemma \ref{lem5} (b), we obtain
\begin{equation}\label{eq24}
	|I_2| \lesssim |x|^{-\hat \beta}\int^{t}_{\frac{t}{2}}\frac{1}{(t-\tau)^{\frac{1}{2}
			(1+\tilde\beta-\hat \beta)}\tau^{\frac{1}{2}(1+\beta-\tilde\beta)}}{\rm d}\tau 
	\simeq |x|^{-\hat \beta}t^{-\frac{1}{2}(\beta-\hat \beta)}.
\end{equation}
From the inequalities \eqref{eq22} and \eqref{eq24}, 
we get $B(u,u)\in K^\beta_{\hat \beta,T}$. From \linebreak $u =e^{t\Delta}u_0+B(u,u)$ 
and Lemma \ref{lem6}, we deduce $u\in K^\beta_{\hat \beta,T}$. This proves the result. Therefore, we get 
$u\in K^\beta_{\hat \beta,T}$ for all $\hat \beta \in [0,\beta]\cap(\tilde \beta-1,\beta]$.\\
We now prove that $u\in K^\beta_{\hat \beta,T}$ for all $\hat \beta \in [0,\beta]$. 
Indeed, in exactly the same way, since $u\in K^\beta_{\hat \beta,T}$ 
for all $\hat \beta \in [0,\beta]\cap(\tilde \beta-1, \beta]$,
it follows that $u\in K^\beta_{\hat \beta,T}$ for all 
$\hat \beta \in [0,\beta]\cap(\tilde \beta-2,\beta]$. By induction, 
we get $u\in K^\beta_{\hat \beta,T}$ for all 
$\hat \beta \in [0,\beta]\cap (\tilde \beta -n, \beta]$ with $n \in \mathbb N$. 
However, there exists a sufficiently large number 
$n$ satisfying $\tilde\beta-n<0$ and therefore $u\in K^\beta_{\hat \beta,T}$ for all 
$\hat \beta \in [0,\beta]$. \linebreak The proof of Lemma \ref{lem11} is complete. \qed
\begin{Lem}\label{lem12}  
	Let $0\leq \beta <d$ be fixed, then for all $\alpha$ and $\tilde \beta$ satisfying
	$$
	\tilde \beta \geq 0, 0 < \alpha <1, \beta -2 < \tilde \beta \leq \beta,
	\  and\ \beta -\tilde \beta -1< \alpha < d-\tilde\beta,
	$$ 
	there exists a positive constant $\delta_{\alpha,\beta,\tilde \beta,d}$ such that for all
	$u_0 \in L^\infty(|x|{\rm d}x) \cap L^\infty(|x|^\beta{\rm d}x)$  
	with ${\rm div}(u_0) = 0$ satisfying
	\begin{equation}\label{eq2}
		\underset{x \in \mathbb R^d, t>0}{\rm sup}\big(|x|^{\alpha}
		t^{\frac{1}{2}(1-  \alpha)} +|x|^{\tilde \beta} 
		t^{\frac{1}{2}(\beta-  \tilde \beta)}\big)|e^{t\Delta}u_0|
		\leq \delta_{\alpha,\beta,\tilde \beta,d},
	\end{equation}
	NSE has a global mild solution $u$ on $(0,\infty) \times \mathbb R^d$ such that
	\begin{equation}\label{eq3}
		\underset{x \in \mathbb R^d, t>0}{\rm sup}\big(|x|+t^{\frac{1}{2}} 
		+|x|^\beta +t^{\frac{\beta}{2}}\big)|u(x,t)| <+\infty.
	\end{equation}
\end{Lem}
\textbf{Proof}. Applying Lemma \ref{lem8}  we deduce that the bilinear operator $B$ is bounded from  $K^1_{\alpha,\infty} \times K^1_{\alpha,\infty}$ 
into $K^1_{\alpha, \infty}$ and from  $K^1_{\alpha,\infty}  \times K^\beta_{\tilde \beta,\infty}$ into $K^\beta_{\tilde \beta,\infty}$. Therefore, 
the  bilinear operator $B$ is bounded from 
$$
(K^1_{\alpha,\infty} \cap K^\beta_{\tilde \beta,\infty})\times (K^1_{\alpha,\infty} \cap K^\beta_{\tilde \beta,\infty}) \ {\rm into}\
(K^1_{\alpha,\infty} \cap K^\beta_{\tilde \beta,\infty}).
$$ 
where the space $K^1_{\alpha,\infty} \cap K^\beta_{\tilde \beta,\infty}$ is equipped with the norm
$$
\big\|u\big\|_{K^1_{\alpha,\infty} \cap K^\beta_{\tilde \beta,\infty}}: = {\rm max}\{\big\|u\big\|_{K^1_{\alpha,\infty}},\big\|u\big\|_{K^\beta_{\tilde \beta,\infty}}\}.
$$
Applying Theorem \ref{th3} to the  bilinear operator $B$, we deduce that there exists a positive constant $\delta_{\alpha,\beta,\tilde \beta,d}$ 
such that for all $u_0 \in L^\infty({|x|{\rm d}x})\cap L^\infty({|x|^\beta{\rm d}x})$  with ${\rm div}(u_0) = 0$ satisfying 
$$
\big\|e^{t\Delta}u_0\big\|_{K^1_{\alpha,\infty} \cap K^\beta_{\tilde \beta,\infty}} \leq \delta_{\alpha,\beta,\tilde \beta,d},
$$
then NSE has a unique mild solution $u$ satisfying
$$
u \in K^1_{\alpha,\infty} \cap K^\beta_{\tilde \beta,\infty}.
$$
Applying Lemmas \ref{lem9}, \ref{lem10}, and \ref{lem11}, we get $u \in K^\beta_{\hat \beta,\infty}$ for all $\hat\beta \in [0,\beta]$ and $u \in K^1_{\alpha,\infty}$ for all $\alpha \in [0,1]$. The proof of Lemma \ref{lem12} is now complete. \qed 
\vskip 0.5cm
{\bf Proof of Theorem \ref{th4}}
\vskip 0.5cm
Since $|x| \leq |x|^\gamma +|x|^\beta$, it follows that 
$$
\big\|f\big\|_{L^\infty}(|x|{\rm d}x) \leq \big\|f\big\|_{L^\infty}(|x|^\gamma{\rm d}x)+\big\|f\big\|_{L^\infty}(|x|^\beta{\rm d}x).
$$
From Lemma \ref{lem2} we have
$$
|x|^\alpha t^{\frac{1}{2}(1-  \alpha)}|e^{t\Delta}u_0| \lesssim \big\|f\big\|_{L^\infty}(|x|{\rm d}x) \leq \big\|f\big\|_{L^\infty}(|x|^\gamma{\rm d}x)+\big\|f\big\|_{L^\infty}(|x|^\beta{\rm d}x),
$$
$$
|x|^{\tilde \gamma}t^{\frac{1}{2}(\gamma-  \tilde \gamma)}|e^{t\Delta}u_0|  \lesssim \big\|f\big\|_{L^\infty}(|x|^\gamma{\rm d}x),\ {\rm and}\ |x|^{\tilde \beta}t^{\frac{1}{2}(\beta-  \tilde \beta)}|e^{t\Delta}u_0| \lesssim  \big\|f\big\|_{L^\infty}(|x|^\beta{\rm d}x).
$$
This proves Theorem \ref{th4}.\qed 

\vskip 0.5cm
{\bf Proof of Theorem \ref{th2}}
\vskip 0.5cm
Since $L^\infty({|x|{\rm d}x})\subset L^\infty({|x|^\gamma{\rm d}x})\cap L^\infty({|x|^\beta{\rm d}x})$, it follows that  $u_0 \in L^\infty({|x|{\rm d}x}) $. Applying Lemma \ref{lem12} then there exists 
a positive constant $\delta_{\alpha,\beta, \tilde \beta,d}$ such that if
$$
\underset{x \in \mathbb R^d, t>0}{\rm sup}\big(|x|^{\alpha}t^{\frac{1}{2}(1-  \alpha)} +|x|^{\tilde \beta} t^{\frac{1}{2}(\beta-  \tilde \beta)}\big)|e^{t\Delta}u_0| \leq \delta_{\alpha,\beta, \tilde \beta,d},
$$
NSE has a global mild solution $u$ on $(0,\infty) \times \mathbb R^d$ such that
$$
\underset{x \in \mathbb R^d, t>0}{\rm sup}\big(|x|+t^{\frac{1}{2}}+|x|^{\beta}+t^{\frac{\beta}{2}}\big)|u(x,t)| <+\infty.
$$
Applying Lemma \ref{lem12} for $\beta =\gamma$ then there exists 
a positive constant $\delta_{\alpha,\gamma, \tilde \gamma, d}$ such that if
$$
\underset{x \in \mathbb R^d, t>0}{\rm sup}\big(|x|^{\alpha}t^{\frac{1}{2}(1-  \alpha)}+|x|^{\tilde \gamma} t^{\frac{1}{2}(\gamma-  \tilde \gamma)}\big)|e^{t\Delta}u_0|\leq \delta_{\alpha,\gamma, \tilde \gamma, d},
$$
NSE has a global mild solution $u$ on $(0,\infty) \times \mathbb R^d$ such that
$$
\underset{x \in \mathbb R^d, t>0}{\rm sup}\big(|x| +t^{\frac{1}{2}}+|x|^{\gamma}+t^{\frac{\gamma}{2}}\big)|u(x,t)| <+\infty.
$$
Therefore, if $u_0$ satisfies the following inequality
$$
\underset{x \in \mathbb R^d, t>0}{\rm sup}\big(|x|^{\tilde \gamma}t^{\frac{1}{2}(\gamma-  \tilde \gamma)} + |x|^\alpha t^{\frac{1}{2}(1-  \alpha)}+|x|^{\tilde\beta} t^{\frac{1}{2}(\beta-  \tilde \beta)}\big)|e^{t\Delta}u_0| \leq
{\min}\{\delta_{\alpha,\beta, \tilde \beta,d},\delta_{\alpha,\gamma, \tilde \gamma, d}\}
$$
NSE has a global mild solution $u$ on $(0,\infty) \times \mathbb R^d$ such that  \eqref{eq5}. \\
The proof of Theorem \ref{th2} is complete.\qed
\vskip 0.5cm
{\bf Proof of Theorem \ref{th2a}}
\vskip 0.2cm 
Applying Lemma \ref{lem8}  we deduce that the bilinear operator $B$ is bounded from  $K^1_{\alpha,\infty} \times K^1_{\alpha,\infty}$ 
into $K^1_{\alpha, \infty}$. Applying Theorem \ref{th3} to the  bilinear operator $B$, we deduce that there exists a positive constant $\delta_{\alpha,d}$ 
such that for all $u_0 \in L^\infty({|x|{\rm d}x})$  with ${\rm div}(u_0) = 0$ satisfying 
$$
\big\|e^{t\Delta}u_0\big\|_{K^1_{\alpha,\infty}} \leq \delta_{\alpha,d},
$$
then NSE has a unique mild solution $u$ satisfying $u \in K^1_{\alpha,\infty}$. Applying \linebreak Lemma \ref{lem9} we have $u \in \underset{\alpha \in [0,1]}{\cap}K^1_{\alpha,\infty}$.\\
We prove that  $u \in K^{\beta}_{\beta, T}$ for all $T \in (0, \infty)$.
Indeed, let $\gamma$ be such that $\gamma \in [1,\beta]\cap(\alpha,\alpha+2)$. Applying Lemma \ref{lem7a} we have $u \in K^{\gamma}_{\alpha,T}$, then using Lemmas \ref{lem10} we get $u \in K^{\gamma}_{\gamma, T}$,  in exactly the same way, using again Lemmas \ref{lem10},  since $u\in K^{\gamma}_{\gamma,T}$ for $\gamma \in [1,\beta]\cap(\alpha,\alpha+2)$, it follows that $u\in K^{\gamma}_{\gamma,T}$ for $\gamma \in [1,\beta]\cap(\alpha,\alpha+4)$. By induction, we get $u\in K^{\gamma}_{\gamma,T}$ for $\gamma \in [1,\beta]\cap(\alpha,\alpha+2n)$. However, there exists a sufficiently large number 
$n$ satisfying $\alpha+2n >\beta$ and therefore  $u \in K^{\beta}_{\beta, T}$. The proof of Theorem  \ref{th2a} is complete. \qed
\vskip 0.5cm
{\bf Acknowledgments}. This research is funded by 
Vietnam National \linebreak Foundation for Science and Technology 
Development (NAFOSTED) under grant number  101.02-2014.50.

\end{document}